\documentclass[10pt,a4paper]{article}

\usepackage{amsmath,amssymb, bbm}
\usepackage{color}
\newcommand{\C}{\mathbb{C}}
\newcommand{\R}{\mathbb{R}}

\def\cA{{\mathcal A}}
\def\cC{{\mathcal C}}

\def\cF{{\mathcal F}}
\def\cP{{\mathcal P}}
\def\cQ{{\mathcal Q}}
\newcommand{\ee}{\varepsilon}

\renewcommand{\div}{{\rm div}\,}

\newcommand{\Sum}{\displaystyle \sum}

\def\d{\partial}

\def\tilde{\widetilde}
\def\hat{\widehat}
\newcommand{\D}{\Delta}
\newcommand{\DF}{\Delta_F}

\newcommand{\n}{\nabla}
\newcommand{\G}{\Gamma}

\newcommand{\Om}{\Omega}
\newcommand{\Ome}{\Omega_\varepsilon}
\newcommand{\tOm}{\tilde{\Omega}_{QG}}

\newcommand{\ve}{v_\ee}
\newcommand{\Ue}{U_\ee}
\newcommand{\Uqg}{U_{\ee,QG}}
\newcommand{\Uosc}{U_{\ee, osc}}
\newcommand{\Uoe}{U_{0,\ee}}
\newcommand{\Uoqg}{U_{0,\ee,QG}}
\newcommand{\Uoosc}{U_{0,\ee, osc}}
\newcommand{\Phie}{\Phi_\ee}
\newcommand{\Thee}{\theta_\ee}

\newcommand{\dOm}{\delta \Omega}

\newtheorem{thm}{Theorem}

\newtheorem{prop}{Proposition}

\newtheorem{rem}{Remark}

\title{Global well-posedness and asymptotics for a penalized Boussinesq-type system without dispersion}

\author{Fr\'ed\'eric Charve\footnote{Universit\'e Paris-Est Cr\'eteil, Laboratoire d'Analyse et de Math\'ematiques Appliqu\'ees (UMR 8050), 61 Avenue du G\'en\'eral de Gaulle, 94 010 Cr\'eteil Cedex (France). E-mail: frederic.charve@u-pec.fr}}

\date{}
\begin{document}

\maketitle

\begin{abstract} J.-Y. Chemin proved the convergence (as the Rossby number $\ee$ goes to zero) of the solutions of the Primitive Equations to the solution of the $3D$ quasi-geostrophic system when the Froude number $F=1$ that is when no dispersive property is available. The result was proved in the particular case where the kinematic viscosity $\nu$ and the thermal diffusivity $\nu'$ are close. In this article we generalize this result for any choice of the viscosities, the key idea is to rely on a special feature of the quasi-geostrophic structure.
\end{abstract}

\section{Introduction}

\subsection{Presentation of the models}

The Primitive Equations we consider in this article (also called Primitive System) are a Boussinesq-type system that describes geophysical flows located in a large scale at the surface of the Earth under the assumption that the vertical motion is much smaller than the horizontal one. Two phenomena have a great influence on geophysical fluids: the rotation of the Earth around its axis and the vertical stratification of the density induced by gravity. The former induces a vertical rigidity in the fluid velocity as described by the Taylor-Proudman theorem, and the latter induces a horizontal rigidity to the fluid density: heavier masses lay under lighter ones.

In order to measure the importance of these two concurrent structures, physicists defined two numbers: the Rossby number $Ro$ and the Froude number $Fr$. We refer to the introduction of \cite{FC, FCpochesLp} for more details and to \cite{BeBo, Cushman, Sadourny, Pedlosky} for an in-depth presentation.

The smaller are these numbers, the more important become these two phenomena and we will consider the Primitive Equations in the whole space, under the Boussinesq approximation and when both phenomena share the same importance i.-e. $Ro=\ee$ and $Fr=\ee F$ with $F\in]0,1]$. In what follows $\ee$ will be called the Rossby number and $F$ the Froude number. The system is then written as follows (we refer to \cite{Chemin2, BMN5} for the model):
\begin{equation}
\begin{cases}
\d_t \Ue +\ve\cdot \n \Ue -L \Ue +\frac{1}{\ee} \cA \Ue=\frac{1}{\ee} (-\n \Phie, 0),\\
\div \ve=0,\\
{\Ue}_{|t=0}=U_{0,\ee}.
\end{cases}
\label{PE}
\tag{$PE_\ee$}
\end{equation}
The unknowns are $\Ue =(\ve, \Thee)=(\ve^1, \ve^2, \ve^3, \Thee)$ (where $\ve$ denotes the velocity of the fluid and $\Thee$ the scalar potential temperature), and $\Phie$ which is called the geopotential. The diffusion operator $L$ is defined by
$$
L\Ue \overset{\mbox{def}}{=} (\nu \D \ve, \nu' \D \Thee),
$$
where $\nu, \nu'>0$ are the kinematic viscosity and the thermal diffusivity. The matrix $\cA$ is defined by
$$
\cA \overset{\mbox{def}}{=}\left(
\begin{array}{llll}
0 & -1 & 0 & 0\\
1 & 0 & 0 & 0\\
0 & 0 & 0 & F^{-1}\\
0 & 0 & -F^{-1} & 0
\end{array}
\right).
$$
We will also precise later the properties satisfied by the sequence of initial data (as $\ee$ goes to zero).

\begin{rem}
 \sl{This system generalises the well-known rotating fluids system and for more precisions we refer to \cite{FCpochesLp}. The fact that $\cA \Ue$ is divided by the Rossby number $\ee$ imposes formal conditions to the limit system as $\ee$ goes to $0$, this term is said to be penalized. The major difference between the classical Navier-Stokes system and \eqref{PE} consists in this penalized term which involves a skew-symmetric matrix, so that for the canonical $\C^4$ inner product and any $L^2$ or $H^s/\dot{H}^s$ inner products, we have $\cA \Ue\cdot \Ue=0$ therefore for all fixed $\ee>0$, any energy method will not "see" these penalized terms and will work as for $(NS)$. Then the Leray and Fujita-Kato theorems are very easily adapted and provide global in time (unique in $2D$) weak solutions if $U_{0,\ee}\in L^2$ and local in time unique strong solutions if $U_{0,\ee} \in \dot{H}^\frac{1}{2}$ (global for small initial data). We refer to Remark \ref{illprepared} for the notion of well/ill-prepared initial data.
}
\label{PenalizeLerayFK}
\end{rem}

\begin{rem}
 \sl{As explained in \cite{FC2, FCpochesLp} two distinct regimes have to be considered regarding the eigenvalues of the linearized system: the case $F\in]0,1[$ where the system features dispersive properties, and the case $F=1$, with simpler operators but where no dispersion occurs. In the dispersive case (see \cite{FC} for weak solutions, \cite{FC2} for strong solutions), using the approach developped by Chemin, Desjardins, Gallagher and Grenier in \cite{CDGG, CDGG2, CDGGbook} for the rotating fluids system, we manage to filter the fast oscillations (going to zero in some norms thanks to Strichartz estimates providing positive powers of the small parameter $\ee$) and prove the convergence to the solution of System \eqref{QG1} below (even for blowing-up ill-prepared initial data as in \cite{FC3, FCpochesLp}, less regular initial data as in \cite{FC4} or with evanescent viscosities as in \cite{FC5}). On the contrary when $F=1$ no dispersion is available and only well-prepared initial data are considered. In addition, in \cite{Chemin2} the asymptotics are obtained only when $\nu$ and $\nu'$ are very close, in \cite{Dragos4} is dealt the inviscid case. We refer also refer to \cite{IMT, KMY, Scro} for results in other context such as periodic domains for example where there is no dispersion, and resonences have to be studied.
 }
\end{rem}

\subsection{The limit system}

We are interested in the asymptotics, as the small parameter $\ee$ goes to zero. Let us recall that in \cite{Chemin2, FC} the limit system, which is a transport-diffusion system coupled with a Biot-Savart inversion law, is first formally obtained, and is called the 3D quasi-geostrophic system:
\begin{equation}
\begin{cases}
\d_t \tOm +\tilde{v}_{QG}.\n \tOm -\G \tOm =0\\
\tilde{U}_{QG}=(\tilde{v}_{QG},\tilde{\theta}_{QG})=(-\partial_2, \partial_1, 0, -F\partial_3) \DF^{-1} \tOm,
\end{cases}
\label{QG1}\tag{QG}
\end{equation}
where the operator $\G$ is defined by:
$$
\G \overset{def}{=} \D \DF^{-1} (\nu \d_1^2 +\nu \d_2^2+ \nu' F^2 \d_3^2),
$$
with $\DF=\d_1^2 +\d_2^2 +F^2 \d_3^2$. Moreover we also have the relation
$$
\tOm=\d_1 \tilde{U}_{QG}^2 -\d_2 \tilde{U}_{QG}^1 -F \d_3 \tilde{U}_{QG}^4 =\d_1 \tilde{v}_{QG}^2 -\d_2 \tilde{v}_{QG}^1 -F \d_3 \tilde{\theta}_{QG}.
$$
\begin{rem}
\sl{The operator $\D_F$ is a simple anisotropic Laplacian but $\G$ is in general a tricky non-local diffusion operator of order 2. In the present article we will focus on the case $F=1$ where $\D_F=\D$ and $\G=\nu \d_1^2 +\nu \d_2^2+ \nu' \d_3^2$. We refer to \cite{FCestimLp, FCpochesLp} for a study of $\G$ in the general case (then neither the Fourier kernel nor the singular integral kernel have a constant sign and no classical result can be used).
}
\end{rem}
\begin{rem}
 \sl{From now on we will consider the very particular case $F=1$. On one hand the operator $\G$ is much simpler, but on the other hand (and as explained for example in \cite{FC, FCestimLp, FCpochesLp}) the system is not dispersive anymore. This lack of dispersive and Strichartz estimates (that were abundantly used in previous works) will force us to use completely different methods, part of them coming from \cite{Chemin2}.
 }
\end{rem}

Led by the limit system we introduce the following decomposition: for any 4-dimensional vector field $U=(v, \theta)$ we define its potential vorticity $\Om(U)$ (here in the case $F=1$):
$$
\Om(U)\overset{def}{=} \d_1 v^2 -\d_2 v^1 -\d_3 \theta,
$$
then its quasi-geostrophic and oscillating (or oscillatory) parts:
\begin{equation}
U_{QG}=\cQ (U) \overset{def}{=} \left(
\begin{array}{r}
-\d_2\\
\d_1\\
0\\
-\d_3
\end{array}
\right) \D^{-1} \Om (U), \quad \mbox{and} \quad U_{osc}=\cP (U) \overset{def}{=} U-U_{QG}.
\end{equation}
As emphasized in \cite{FC,FC5} this is an orthogonal decomposition of 4-dimensional vector fields (similar to the Leray orthogonal decomposition into divergence-free and gradient vector fields) and if $\cQ$ and $\cP$ are the associated orthogonal projectors on the quasi-geostrophic or oscillating fields, they satisfy (see \cite{Chemin2, FC, FC2}):
\begin{prop}
\sl{With the same notations, for any function $U=(v, \theta)$ we have:
\begin{enumerate}
\item $\cP$ and $\cQ$ are pseudo-differential operators of order 0.
\item For any $s\in\R$, $(\cP(U)|\cQ(U))_{\dot{H}^s} =(\cA U|\cP(U))_{\dot{H}^s}=0$.
\item The same is true for nonhomogeneous Sobolev spaces.
\item $\cP(U)=U \Longleftrightarrow \cQ(U)=0\Longleftrightarrow \Om(U)=0$.
\item $\cQ(U)=U \Longleftrightarrow \cP(U)=0\Longleftrightarrow$ there exists a scalar function $\Phi$ such that $U=(-\d_2,\d_1,0,-\d_3) \Phi$. Such a vector field is said to be quasi-geostrophic and is divergence-free.
\item If $U=(v, \theta)$ is a quasi-geostrophic vector field, then $v\cdot \n \Om(U)=\Om(v\cdot \n U)$.
\item If $U$ is a quasi-geostrophic vector field, then $\G U=\cQ (L U)$.
\end{enumerate}
\label{propdecomposcqg}
}
\end{prop}
Thanks to this, System \eqref{QG1} can for example be rewritten into the following velocity formulation:
\begin{equation}
\begin{cases}
\d_t \tilde{U}_{QG} +\tilde{v}_{QG}.\n \tilde{U}_{QG} -L\tilde{U}_{QG}= \cP \tilde{\Phi}_{QG},\\
\tilde{U}_{QG}=\mathcal{Q} (\tilde{U}_{QG}), \mbox{ (or equivalently } \mathcal{P} (\tilde{U}_{QG})=0),\\
\tilde{U}_{QG}|_{t=0}=\tilde{U}_{0,QG}.
\label{QG2}\tag{$QG_2$}
\end{cases}
\end{equation}
Back to System \eqref{PE}, if we introduce $\Ome=\Om(\Ue)$, $\Uqg=\cQ(\Ue)$ and $\Uosc=\cP(\Ue)$, they satisfy the following systems (see \cite{FC} for details):
\begin{equation}
\d_t \Ome +\ve\cdot \n \Ome -\G \Ome= (\nu-\nu')\D \d_3 \theta_{\ee, osc}+q_\ee,
\label{systomega}
\end{equation}
where $q_\ee$ is defined by
\begin{multline}
q_\ee =q(\Uosc, \Ue) =\d_3 v_{\ee, osc}^3(\d_1 \ve^2-\d_2 \ve^1) -\d_1 v_{\ee, osc}^3 \d_3 \ve^2 +\d_2 v_{\ee, osc}^3 \d_3 \ve^1\\
+\d_3 v_{\ee, QG}\cdot \n \theta_{\ee, osc} +\d_3 v_{\ee, osc}\cdot \n \Thee,
\end{multline}
and
\begin{multline}
\d_t \Uosc-(L-\frac{1}{\ee}\mathbb{P} \cA) \Uosc = -\mathbb{P}(\ve\cdot \n \Ue) -\left(
\begin{array}{c}
-\d_2\\
\d_1\\
0\\
-\d_3
\end{array}
\right) \D^{-1} \bigg( -\ve\cdot \n \Ome +q_\ee\bigg)\\
+(\nu-\nu') \d_3\left(
\begin{array}{c}
\d_2 \theta_\ee\\
-\d_1 \theta_\ee\\
0\\
\d_1 v_\ee^2-\d_2 v_\ee^1
\end{array}
\right).
\label{systosc}
\end{multline}

\begin{rem}
\sl{For more conciseness and without any loss of generality, we will write in what follows
$$
q_\ee=\n \Uosc \cdot \n \Ue.
$$
}
\end{rem}

\begin{rem}
\sl{It is natural to investigate the link between the quasi-geostrophic/oscillating parts decomposition of the initial data and the asymptotics when $\ee$ goes to zero. This leads to the notion of well-prepared/ill-prepared initial data depending on the fact that the initial data is already close or not to the quasi-geostrophic structure, i.-e. when the initial oscillating part is small/large (or going to zero/blowing up as $\ee$ goes to zero). We refer to \cite{FCpochesLp} for more details about this subject. For example in \cite{FC3, FCestimLp, FCpochesLp} we focussed on the case $F\in]0,1[$ for very ill-prepared cases in the sense that the initial oscillating part norm goes to infinity as $\ee$ goes to zero, a way to balance these large norms was to take advantage of the dispersive estimates satisfied by the oscillating part, providing positive powers of $\ee$. On the contrary when $F=1$, as in \cite{Chemin2}, we will consider well-prepared initial data.}
\label{illprepared}
\end{rem}

\subsection{Statement of the main results}

The aim of the present article is to generalize the results of Chemin from \cite{Chemin2}, which were obtained only in the case where $\nu \sim \nu'$. To the best of our knowledge, this study has never been investigated any further in the non-dispersive case $F=1$. First let us define the family of spaces $\dot{E}_T^s$ for $s\in \R$,
$$
\dot{E}_T^s=\mathcal{C}_T(\dot{H}^s) \cap L_T^2(\dot{H}^{s+1}),
$$
endowed with the following norm (see the appendix for notations):
$$
\|f\|_{\dot{E}_T^s}^2 \overset{def}{=}\|f\|_{L_T^\infty \dot{H}^s}^2+\min(\nu, \nu')\int_0^T \|f(\tau)\|_{\dot{H}^{s+1}}^2 d\tau.
$$
When $T=\infty$ we denote $\dot{E}^s$ and the corresponding norm is over $\R_+$ in time.\\

\begin{thm} (Global existence and uniqueness)
 \sl{There exists a positive constant $\mathbb{C}$ such that for any initial data $U_0\in H^1=L^2 \cap \dot{H}^1$, if
 \begin{equation}
  \begin{cases}
  \vspace{0.2cm}
  \displaystyle{\|U_{0,osc}\|_{\dot{H}^{-1}} \leq \frac{1}{\mathbb{C}^2}\frac{\min(\nu, \nu')^4}{\|U_0\|_{\dot{H}^1}^3} \exp\left(-\mathbb{C} \frac{\|U_0\|_{L^2} \|U_0\|_{\dot{H}^1}}{\min(\nu, \nu')^2}\right),}\\
  \displaystyle{\ee \leq \frac{1}{\mathbb{C}^2}\frac{\min(\nu, \nu')^4}{\|U_0\|_{\dot{H}^1}^4 \Big(\|U_0\|_{\dot{H}^\frac{1}{2}}+\max(\nu,\nu')\Big)} \exp\left(-\mathbb{C}\frac{ \|U_0\|_{L^2} \|U_0\|_{\dot{H}^1}}{\min(\nu, \nu')^2}\right),}
 \end{cases}
 \label{Condepsosc}
 \end{equation}
then System \eqref{PE} has a unique global solution in the space $\dot{E}^1$ and we have $\|\Ue\|_{\dot{E}^1}\leq 2\|U_0\|_{\dot{H}^1}$.
}
 \label{Th1}
\end{thm}

\begin{rem}
 \sl{Of course, due to the Leray estimates, we obtain that in fact all norms in $\dot{E}^s$ for $s\in[0,1]$ are uniformly bounded and more precisely, for all $t\geq 0$,
 $$
 \|\Ue(t)\|_{\dot{H}^s}^2+\min(\nu, \nu')\int_0^t \|\Ue(\tau)\|_{\dot{H}^{s+1}}^2 d\tau \leq \|U_0\|_{L^2}^{2(1-s)}(2\|U_0\|_{\dot{H}^1})^{2s}.
 $$
 }
\end{rem}
The second result is devoted to the asymptotics as $\ee\rightarrow 0$:
\begin{thm} (Convergence)
 \sl{Let $(\Uoe)_{\ee \in]0,\ee_0]}$ be a family of initial data, uniformly bounded in $L^2\cap \dot{H}^1$. Assume in addition that $\Uoosc \in \dot{H}^{-1}$ and there exists a quasi-geostrophic function $\tilde{U}_{0,QG}$ and $\delta>0$ such that:
 $$
 \begin{cases}
 \vspace{0.2cm}
  \Uoosc \underset{\ee \rightarrow 0}{\longrightarrow} 0 \mbox{ in } \dot{H}^{-1},\\
  \vspace{0.2cm}
  \Uoosc \mbox{ is uniformly bounded in }\dot{H}^{1+\delta},\\
  \Uoqg \underset{\ee \rightarrow 0}{\longrightarrow} \tilde{U}_{0,QG} \mbox{ in } \dot{H}^1,
 \end{cases}
 $$
 then there exists $\ee_1\in]0,\ee_0]$ such that for any $\ee<\ee_1$, the assumptions from the previous theorem are fulfilled so that System \eqref{PE} admits a unique global solution $\Ue$ and the family $(\Ue)_{\ee \in]0,\ee_1]}$ converges to the unique global solution $\tilde{U}_{QG}$ of System \eqref{QG1} with initial data $\tilde{U}_{0,QG}$ in the following sense:
 $$
 \begin{cases}
  \vspace{0.2cm}
    \Uosc \underset{\ee \rightarrow 0}{\longrightarrow} 0 \mbox{ in } \dot{E}^s\mbox{ for any } s\in [-1,1[,\\
    \Uqg-\tilde{U}_{QG} \underset{\ee \rightarrow 0}{\longrightarrow} 0 \mbox{ in } \dot{E}^s\mbox{ for any } s\in ]0,1].\\
 \end{cases}
 $$
 }
 \label{Th2}
\end{thm}

\begin{rem}
 \sl{The extra regularity assumption $1+\delta$ on the oscillating part is needed to prove the convergence of the quasi-geotrophic part.}
\end{rem}

\section{Proof of Theorem \ref{Th1}}

In this section, we will follow some parts of the method from \cite{Chemin2}, we will also explain which parts are useless if we do not assume $\nu\sim \nu'$ anymore, and explain how (as in \cite{FCpochesLp}) the quasi-geostrophic structure can once more help obtaining the result.
\\

Let us recall that according to Remark \ref{PenalizeLerayFK}, as $U_0\in L^2$ we have a global weak Leray-type solution $\Ue \in L^\infty (L^2) \cap L^2(\dot{H}^1)$. Moreover as in addition $U_0\in \dot{H}^1$ then it also belongs to $\dot{H}^\frac{1}{2}$ and thanks to the Fujita and Kato theorem, and the weak-strong uniqueness, $\Ue$ is also the unique local strong solution in $\dot{E}_T^\frac{1}{2}$ for all $T <T_\ee^*$ which is the maximal lifespan.

Moreover we recall the classical propagation of regularity estimates (we refer for example to \cite{Chemin1, Dbook}, or \cite{Chemin2} section $2$, for more details in the Navier-Stokes setting, and the proofs which only rely on energy estimates, are still valid in our case), for all $s\in ]-1,1]$, for all $t$,
\begin{equation}
 \|\Ue(t)\|_{\dot{H}^s}^2 +\min(\nu, \nu') \int_0^t \|\n\Ue(\tau)\|_{\dot{H}^s}^2 d\tau \leq \|U_0\|_{\dot{H}^s}^2 \exp \left(\frac{C}{\min(\nu, \nu')}\int_0^t \|\Ue(\tau)\|_{\dot{H}^\frac{3}{2}}^2 d\tau\right),
 \label{estimapriori}
\end{equation}
Thanks to this we also have $\Ue \in \dot{E}_T^1$ for all $T<T_\ee^*$ and finally we have the classical blowup criterion:
\begin{equation}
T_\ee^*<\infty \Rightarrow \int_0^{T_\ee^*} \|\nabla \Ue(\tau)\|_{\dot{H}^\frac{1}{2}}^2 d\tau =\infty.
 \label{Blowupcrit}
\end{equation}

\begin{rem}
 \sl{Let us emphasize that the previous blowup criterion is given by results adapted to the general Navier-Stokes system (the penalization term is invisible to them so we can use them) but according to our computations, and as also observed in \cite{Chemin2}, it is in fact the smaller quantity $\int_0^{T_\ee^*} \|\nabla \Uosc(\tau)\|_{\dot{H}^\frac{1}{2}}^2 d\tau$ that controls the lifespan (see \eqref{controleoscb}).}
 \label{controleosc}
\end{rem}

\subsection{First step: estimates on $\d_t \Ue+\frac{1}{\ee}\cA \Uosc$} 
The first important idea in \cite{Chemin2} consists in changing the formulation of System \eqref{PE}: taking the divergence of the first three equations from \eqref{PE} we get:
$$
-\frac{1}{\ee} \D \Phie =\div(\ve\cdot \n \Ue)+\frac{1}{\ee} \div(\cA \Ue) =\Sum_{1\leq i,j\leq 3} \d_i \d_j(\ve^i \ve^j) -\frac{1}{\ee} \Ome.
$$
Observing that $\cA \Uqg=-(\n \D^{-1} \Ome, 0)$, we end up with the following equivalent formulation of \eqref{PE}:
\begin{equation}
\begin{cases}
\d_t \Ue +\ve\cdot \n \Ue -L \Ue +\frac{1}{\ee} \cA \Uosc=(\n p_{\ee, osc},0) \overset{def}{=} \Big(\n \Sum_{1\leq i,j\leq 3} \d_i \d_j \D^{-1}(\ve^i \ve^j), 0\Big),\\
\div \ve=0,\\
{\Ue}_{|t=0}=U_{0,\ee}.
\end{cases}
\label{PE2}
\tag{$PE_{\ee,2}$}
\end{equation}
From this reformulation, Chemin connects the $\dot{H}^{-1}$ or $L^2$ norm of the block $\ee^{-1} \Uosc$ to $\d_t \Ue$ through the following proposition, that we state here in the general setting for $\nu, \nu'$:
\begin{prop}
 \sl{There exists a constant $C>0$ such that for any $\Ue$ solution of \eqref{PE} in $\cC_T(L^2\cap \dot{H}^1)$, the following estimates hold for all $t\leq T$:
 \begin{equation}
  \begin{cases}
  \vspace{0.2cm}
   \displaystyle{\Big\|\d_t \Ue(t)+ \frac{1}{\ee}\cA \Uosc(t)\Big\|_{\dot{H}^{-1}} \leq \|\Ue(t)\|_{\dot{H}^1} \big(\max(\nu,\nu') +C\|\Ue(t)\|_{\dot{H}^\frac{1}{2}}\big),}\\
   \displaystyle{\Big\|\d_t \Ue(t)+ \frac{1}{\ee}\cA \Uosc(t)\Big\|_{L^2} \leq \|\Ue(t)\|_{\dot{H}^2}\big(\max(\nu,\nu') +C\|\Ue(t)\|_{\dot{H}^\frac{1}{2}}}\big).
  \end{cases}
  \label{UoscdtU}
 \end{equation}
 }
\end{prop}
Then we derive in time System \eqref{PE2} and obtain (we easily adapt \cite{Chemin2}, Section 2 step 3, and skip details) that for all $t$:
\begin{equation}
 \|\d_t \Ue(t)\|_{\dot{H}^{-1}}^2 +\min(\nu, \nu') \int_0^t \|\d_t \Ue(\tau)\|_{L^2}^2 d\tau \leq \|\d_t \Ue(0)\|_{\dot{H}^{-1}}^2 \exp \left(\frac{C}{\min(\nu, \nu')}\int_0^t \|\Ue(\tau)\|_{\dot{H}^\frac{3}{2}}^2 d\tau\right).
\end{equation}
Using \eqref{UoscdtU} to roughly estimate $\|\d_t \Ue(0)\|_{\dot{H}^{-1}}$ we obtain:
\begin{multline}
  \|\d_t \Ue(t)\|_{\dot{H}^{-1}}^2 +\min(\nu, \nu') \int_0^t \|\d_t \Ue(\tau)\|_{L^2}^2 d\tau\\
 \leq \Big(\frac{1}{\ee}\|U_{0,osc}\|_{\dot{H}^{-1}} +\|U_0\|_{\dot{H}^1} \big(\max(\nu,\nu') +C\|U_0\|_{\dot{H}^\frac{1}{2}}\big)\Big)^2 \exp \left(\frac{C}{\min(\nu, \nu')}\int_0^t \|\Ue(\tau)\|_{\dot{H}^\frac{3}{2}}^2 d\tau\right).
\end{multline}
Then, writing
$$
\|\frac{1}{\ee}\cA \Uosc(t)\|_{\dot{H}^{-1}} \leq  \|\frac{1}{\ee}\cA \Uosc(t) +\d_t \Ue(t)\|_{\dot{H}^{-1}} +\|\d_t \Ue(t)\|_{\dot{H}^{-1}},
$$
and thanks once more to \eqref{UoscdtU}, we get:
\begin{multline}
\|\Uosc(t)\|_{\dot{H}^{-1}} \leq \ee \|\Ue(t)\|_{\dot{H}^1} \big(\max(\nu,\nu') +C\|\Ue(t)\|_{\dot{H}^\frac{1}{2}}\big)\\
+\Big(\|U_{0,osc}\|_{\dot{H}^{-1}} +\ee\|U_0\|_{\dot{H}^1} \big(\max(\nu,\nu') +C\|U_0\|_{\dot{H}^\frac{1}{2}}\big)\Big) \exp \left(\frac{C}{2\min(\nu, \nu')}\int_0^t \|\Ue(\tau)\|_{\dot{H}^\frac{3}{2}}^2 d\tau\right).
\end{multline}
Thanks to \eqref{estimapriori} we end up with:
\begin{multline}
\|\Uosc(t)\|_{\dot{H}^{-1}} \leq \Big(\|U_{0,osc}\|_{\dot{H}^{-1}} +2\ee\|U_0\|_{\dot{H}^1} \big(\max(\nu,\nu') +C\|U_0\|_{\dot{H}^\frac{1}{2}}\big)\Big)\\
\times \exp \left(\frac{C}{\min(\nu, \nu')}\int_0^t \|\Ue(\tau)\|_{\dot{H}^\frac{3}{2}}^2 d\tau\right).
\label{estimosc-1}
\end{multline}
Similarly, we obtain that for all $t$:
\begin{multline}
\min(\nu, \nu') \int_0^t \|\Uosc(\tau)\|_{L^2}^2 d\tau \leq 4\Big(\|U_{0,osc}\|_{\dot{H}^{-1}} +\ee\|U_0\|_{\dot{H}^1} \big(\max(\nu,\nu') +C\|U_0\|_{\dot{H}^\frac{1}{2}}\big)\Big)^2\\
\times \exp \left(\frac{C}{\min(\nu, \nu')}\int_0^t \|\Ue(\tau)\|_{\dot{H}^\frac{3}{2}}^2 d\tau\right).
\label{estimosc0}
\end{multline}

\subsection{Second step: Energy estimates}

In \cite{Chemin2}, Chemin uses the quasi-geostrophic/oscillating orthogonal decomposition of the solution and estimates each part by energy methods applied to Systems \eqref{systomega} and \eqref{systosc}, combining them through:
$$
\|\Ue\|_{\dot{H}^1}^2 \sim \|\Uqg\|_{\dot{H}^1}^2 +\|\Uosc\|_{\dot{H}^1}^2 \sim \|\Ome\|_{L^2}^2 +\|\Uosc\|_{\dot{H}^1}^2 ,
$$
More precisely, instead of \eqref{systosc}, using the changes in the pressure term that lead to System \eqref{PE2}, Chemin studied the following system (that we write here in a more accurate QG/osc decomposition):
\begin{multline}
 \d_t \Uosc +\ve \cdot \n \Uosc + \frac{1}{\ee} \cA \Uosc=(\n p_{\ee, osc}, 0) +\left(
 \begin{array}{c}
  [\ve\cdot \n, \d_2 \D^{-1}]\\ -[\ve\cdot \n, \d_1 \D^{-1}]\\ 0\\ -[\ve\cdot \n, \d_3 \D^{-1}]
 \end{array}
\right)\Ome\\
 +\left(
 \begin{array}{c}
  -\d_2\\ \d_1\\ 0\\ -\d_3
 \end{array}
\right) \Big(-(\nu-\nu')\d_3 \Thee +\D^{-1} q_\ee\Big) +(\nu-\nu')\left(
\begin{array}{c}
 0\\ 0\\ 0\\ \d_3
\end{array}
\right) \Ome.
\label{systosc2}
\end{multline}
In the general case where $\nu-\nu'$ is not assumed to be small, this method is bound to fail. Indeed taking the $L^2$-inner product of \eqref{systomega} with $\Ome$, we obtain:
$$
\frac{1}{2}\frac{d}{dt}\|\Ome\|_{L^2} +\min(\nu, \nu') \|\n \Ome\|_{L^2} \leq (\nu-\nu')(\partial_3 \D \Theta_{\ee, osc},\Ome)_{L^2} + (q_\ee,\Ome)_{L^2}.
$$
As emphasized in \cite{FCpochesLp} the first term of the right-hand side features too many derivatives for us to be able to prove it is small, the best we can hope for is to bound it by $\|\Uosc\|_{\dot{H}^2}\|\Ome\|_{\dot{H}^1}$ and obtain that this term is at most bounded whereas we need it to go to zero with $\ee$.\\

Any hope to neutralize this term thanks to the $\dot{H}^1$-innerproduct of \eqref{systosc2} by $\Uosc$ is also completely out of reach as the last term in \eqref{systosc2} will produce the exact same term in the right-hand side of the energy estimate instead of a cancellation.

The only way to get rid of this difficulty is to go back to the original system and use \emph{another very important feature of the quasi-geostrophic decomposition}. As was first obtained in \cite{FCpochesLp} (in the case $F\neq 1$, see Section 4, let us emphasize that even in this case, dispersion cannot help when it is about estimating in $L^2$), we can take advantage of Point 6 from Proposition \ref{propdecomposcqg} (this property was first observed in \cite{FC2}) and gain important cancellations. Taking the $\dot{H}^1$-inner product of \eqref{PE} with $\Ue$, we obtain that
\begin{equation}
 \frac{1}{2} \frac{d}{dt} \|\Ue\|_{\dot{H}^1}^2 -(L \Ue| \Ue)_{\dot{H}^1}=-(v_\ee\cdot \n\Ue| \Ue)_{\dot{H}^1}.
\label{PSUe}
 \end{equation}
As before $-(L \Ue| \Ue)_{\dot{H}^1} \geq \min(\nu, \nu') \|\n \Ue\|_{\dot{H}^1}^2$ and thanks to the decomposition $\Ue=\Uosc+\Uqg$ we develop the right-hand side as follows:
\begin{multline}
(v_\ee\cdot \n\Ue| \Ue)_{\dot{H}^1} =(v_\ee\cdot \n\Ue| \Uosc)_{\dot{H}^1} +(v_\ee\cdot \n\Uosc| \Uqg)_{\dot{H}^1} +(v_{\ee,osc}\cdot \n\Uqg| \Uqg)_{\dot{H}^1}\\
+(v_{\ee,QG}\cdot \n\Uqg| \Uqg)_{\dot{H}^1} \overset{def}{=} A_1+ A_2+A_3+A_4.
\end{multline}
Then, as emphasized in \cite{FCpochesLp}, the last term (which is the most dangerous term as it does not involve any occurence of the evanescent oscillating part, and therefore has no reason, at first sight, to go to zero as $\ee$ goes to zero) is \emph{equal to zero}, this is the key property allowing us to complete the proof. To show this we simply use the following elementary computation related to the quasi-geostrophic decomposition: for any function $f$, we have (here in the particular case $F=1$):
$$
(f|\Uqg)_{\dot{H}^1} =-(f|\left(
\begin{array}{c}
-\d_2\\
\d_1\\
0\\
-\d_3
\end{array}
\right)\Ome)_{L^2} =(\Om(f)|\Ome)_{L^2}.
$$
Then, thanks to Point $6$ from Proposition \ref{propdecomposcqg} and the fact that $\div v_{\ee,QG}=0$, we obtain:
$$
(v_{\ee,QG}\cdot \n\Uqg| \Uqg)_{\dot{H}^1}=(\Om(v_{\ee,QG}\cdot \n\Uqg)| \Ome)_{L^2} =(v_{\ee,QG}\cdot \n\Ome| \Ome)_{L^2}=0.
$$
Next let us estimate the first three terms. Thanks to the Sobolev injections or product laws (see appendix),
\begin{equation}
\begin{cases}
\vspace{0.2cm}
|A_1|\leq \|\ve\cdot \n \Ue\|_{\dot{H}^\frac{1}{2}} \|\Uosc\|_{\dot{H}^\frac{3}{2}} \leq  \|\ve\|_{\dot{H}^1} \cdot \|\Ue\|_{\dot{H}^2} \cdot \|\Uosc\|_{\dot{H}^\frac{3}{2}},\\
|A_3|\leq \|\ve\cdot \n \Uosc\|_{L^2} \|\Uqg\|_{\dot{H}^2} \leq \|\ve\|_{\dot{H}^1} \cdot \|\n \Uosc\|_{\dot{H}^\frac{1}{2}} \cdot \|\Ue\|_{\dot{H}^2}.
\end{cases}
\end{equation}
Next using $\div{v_{\ee, osc}}=0$ for $A_2$,
\begin{multline}
|A_2|=|\Sum_{i=1,...,3} \int_{\R^3} \d_i (v_{\ee, osc}\cdot \n \Uqg) \cdot \d_i \Uqg dx| =|\Sum_{i=1,...,3} \int_{\R^3} \d_i v_{\ee, osc}\cdot \n \Uqg \cdot \d_i \Uqg dx|\\
\leq \|\n \Uosc\|_{L^3} \cdot \|\n\Uqg\|_{L^6} \cdot \|\n \Uqg\|_{L^2} \leq \|\n \Uosc\|_{\dot{H}^\frac{1}{2}} \cdot \|\Ue\|_{\dot{H}^2} \cdot \|\Ue\|_{\dot{H}^1}.
\end{multline}
Plugging this into \eqref{PSUe} and thanks to the classical estimates $2ab\leq a^2+b^2$ we get:
$$
\frac{d}{dt} \|\Ue\|_{\dot{H}^1}^2 +\min(\nu, \nu')\|\n \Ue\|_{\dot{H}^1}^2 \leq \frac{C}{\min(\nu, \nu')} \|\Ue\|_{\dot{H}^1}^2 \cdot \|\Uosc\|_{\dot{H}^\frac{3}{2}}^2,
$$
and thanks to the Gronwall estimates, we end up (as announced in Remark \ref{controleosc}) with for all $t<T_\ee^*$: 
\begin{equation}
\|\Ue(t)\|_{\dot{H}^1}^2 +\min(\nu, \nu') \int_0^t\|\n \Ue(\tau)\|_{\dot{H}^1}^2 d\tau \leq \|\Uoe\|_{\dot{H}^1}^2 \exp \Big(\frac{C}{\min(\nu, \nu')} \int_0^t \|\Uosc(\tau)\|_{\dot{H}^\frac{3}{2}}^2d\tau\Big).
\label{controleoscb}
\end{equation}
Consequently, if
$$
\int_0^t \|\Uosc(\tau)\|_{\dot{H}^\frac{3}{2}}^2d\tau \leq \frac{\ln 2}{C} \min(\nu, \nu'),
$$
then we have
\begin{equation}
\|\Ue(t)\|_{\dot{H}^1}^2 +\min(\nu, \nu') \int_0^t\|\n \Ue(\tau)\|_{\dot{H}^1}^2 d\tau \leq 2\|\Uoe\|_{\dot{H}^1}^2,
\end{equation}
and thanks to the Leray estimates together with interpolation:
$$
\int_0^t\|\n \Ue(\tau)\|_{\dot{H}^\frac{3}{2}}^2 d\tau \leq \Big( \int_0^t\|\n \Ue(\tau)\|_{\dot{H}^1}^2 d\tau \Big)^\frac{1}{2} \Big( \int_0^t\|\n \Ue(\tau)\|_{\dot{H}^2}^2 d\tau \Big)^\frac{1}{2} \leq \sqrt{2} \frac{\|U_0\|_{L^2} \|U_0\|_{\dot{H}^1}}{\min(\nu, \nu')},
$$
which allows us to state some analoguous to Corollary $2.1$ from \cite{Chemin2}:
\begin{prop}
 \sl{Let $\Ue\in \cC([0,T_\ee^*[, \dot{H}^1) \cap L_{loc}^2([0,T_\ee^*[, \dot{H}^2)$ be a solution of \eqref{PE}. If there exists some $T_\ee>0$ such that
$$
\int_0^{T_\ee} \|\Uosc(\tau)\|_{\dot{H}^\frac{3}{2}}^2d\tau \leq \frac{\ln 2}{C} \min(\nu, \nu'),
$$
then for all $t\leq T_\ee$,
\begin{equation}
 \begin{cases}
\vspace{0.2cm}
 \displaystyle{\|\Ue(t)\|_{\dot{H}^1}^2 +\min(\nu, \nu') \int_0^t\|\Ue(\tau)\|_{\dot{H}^2}^2 d\tau \leq 2\|U_0\|_{\dot{H}^1}^2,}\\
 \displaystyle{\int_0^t\|\n \Ue(\tau)\|_{\dot{H}^\frac{3}{2}}^2 d\tau \leq \sqrt{2}\frac{\|U_0\|_{L^2} \|U_0\|_{\dot{H}^1}}{\min(\nu, \nu')}.}
\end{cases}
\label{estimboot}
\end{equation}
 }
 \label{propuosc}
\end{prop}

\subsection{Third step: boostrap and proof of Theorem \ref{Th1}}

From the previous estimates we will develop a boostrap argument to prove the first theorem. For $\ee>0$ fixed (and which will be precised later), we consider the unique local strong solution $\Ue$ built in the beginning of the previous section. We recall that $\Ue \in \dot{E}_t^s$ for all $t<T_\ee^*$ and $s\in [0,1]$ and that in addition, thanks to the Leray estimates, $\|\Ue\|_{\dot{E}_t^0} \leq \|U_{0}\|_{L^2}$ for all $t<T_\ee^*$.

Let us define:
\begin{equation}
 T_\ee =\sup \{t\in ]0, T_\ee^*[, \int_0^t \|\Uosc(\tau)\|_{\dot{H}^\frac{3}{2}}^2d\tau \leq \frac{\ln 2}{C} \min(\nu, \nu') \}.
\end{equation}
Obviously, $T_\ee\in ]0, T_\ee^*]$, and if $T_\ee= T_\ee^*$ then by Proposition \ref{propuosc} and the blowup criterion from \eqref{Blowupcrit} we immediately obtain that $\Ue$ is global and \eqref{estimboot} becomes valid for any time.

If not then, as $T_\ee\in ]0, T_\ee^*[$, it is finite and
\begin{equation}
 \int_0^{T_\ee} \|\Uosc(\tau)\|_{\dot{H}^\frac{3}{2}}^2d\tau =\frac{\ln 2}{C} \min(\nu, \nu').
\label{Boot}
\end{equation}
Moreover for all $t\leq T_\ee$, combining \eqref{estimosc-1}, \eqref{estimosc0} with \eqref{estimboot}, we get that for all $t\leq T_\ee$
\begin{multline}
\max\left(\|\Uosc(t)\|_{\dot{H}^{-1}}^2, \min(\nu, \nu') \int_0^t \|\Uosc(\tau)\|_{L^2}^2 d\tau \right)\\
\leq 4\Big(\|U_{0,osc}\|_{\dot{H}^{-1}} +\ee\|U_0\|_{\dot{H}^1} \big(\max(\nu,\nu') +C\|U_0\|_{\dot{H}^\frac{1}{2}}\big)\Big)^2 \exp \left(C \frac{\|U_0\|_{L^2} \|U_0\|_{\dot{H}^1}}{\min(\nu, \nu')^2} \right).
\label{estimosc0b}
\end{multline}
As $\Uosc=\cP \Ue$ (with $\cP$ homogeneous Fourier multiplier of order zero), the first majoration from \eqref{estimboot} is also true for $\Uosc$, and combining it with the previous estimates through interpolation ($3/2=(1-\theta)\cdot 0+\theta \cdot 2$ with $\theta=3/4$) we obtain:
\begin{multline}
 \int_0^t \|\Uosc(\tau)\|_{\dot{H}^\frac{3}{2}}^2d\tau \leq \frac{2^\frac{5}{4}\|U_0\|_{\dot{H}^1}^\frac{3}{2}}{\min(\nu, \nu')} \Big(\|U_{0,osc}\|_{\dot{H}^{-1}} +\ee\|U_0\|_{\dot{H}^1} \big(\max(\nu,\nu') +C\|U_0\|_{\dot{H}^\frac{1}{2}}\big)\Big)^\frac{1}{2}\\
 \times \exp \left(C \frac{\|U_0\|_{L^2} \|U_0\|_{\dot{H}^1}}{\min(\nu, \nu')^2} \right).
\end{multline}
This quantity is therefore less than $\frac{\ln 2}{2C} \min(\nu, \nu')$ as soon as
$$
\|U_{0,osc}\|_{\dot{H}^{-1}} +\ee\|U_0\|_{\dot{H}^1} \big(\max(\nu,\nu') +C\|U_0\|_{\dot{H}^\frac{1}{2}}\big) \leq \frac{\ln 2)^2}{C^2 2^\frac{9}{2}} \frac{\min(\nu, \nu')^4}{\|U_0\|_{\dot{H}^1}^3}\exp \left(\frac{C}{\sqrt{2}} \frac{\|U_0\|_{L^2} \|U_0\|_{\dot{H}^1}}{\min(\nu, \nu')^2} \right),
$$
which is realized when $\ee$ and $U_{0,osc}$ satisfy \eqref{Condepsosc} and then we have proved that for any $t\leq T_\ee$, we have in fact
$$
\int_0^t \|\Uosc(\tau)\|_{\dot{H}^\frac{3}{2}}^2d\tau \leq \frac{\ln 2}{2C} \min(\nu, \nu'),
$$
which contradicts \eqref{Boot} and the definition of $T_\ee$, then $T_\ee=T_\ee^*=\infty$ and \eqref{estimboot} and \eqref{estimosc0b} are valid for any time and Theorem \ref{Th1} is proved $\blacksquare$.

By interpolation we deduce that for all $s\in[-1,1[$,
\begin{multline}
 \|\Uosc\|_{\dot{E}^s} \leq C\Big(\|U_{0,osc}\|_{\dot{H}^{-1}} +\ee\|U_0\|_{\dot{H}^1} \big(\max(\nu,\nu') +C\|U_0\|_{\dot{H}^\frac{1}{2}}\big)\Big)^\frac{1-s}{2}\\
 \times \exp \left(C \frac{\|U_0\|_{L^2} \|U_0\|_{\dot{H}^1}}{\min(\nu, \nu')^2} \right) \|U_0\|_{\dot{H}^1}^\frac{1+s}{2}.
\end{multline}
As $\|U_{0,osc}\|_{\dot{H}^{-1}}$ is constant, this estimate is useless, but if we consider the initial data from Theorem \ref{Th2}, if $\|\Uoe\|_{H^1} \leq C_0$, then 
$$
\frac{1}{\mathbb{C}^2}\frac{\min(\nu, \nu')^4}{\|U_0\|_{\dot{H}^1}^3} \exp\left(-\mathbb{C} \frac{\|U_0\|_{L^2} \|U_0\|_{\dot{H}^1}}{\min(\nu, \nu')^2}\right) \geq \frac{\min(\nu, \nu')^4}{\mathbb{C}^2 C_0^3} \exp\left(-\frac{\mathbb{C}C_0^2}{\min(\nu, \nu')^2}\right)>0,
$$
which bounds $\|\Uoosc\|_{\dot{H}^{-1}}$ when $\ee$ is small enough, thanks to the assumptions. Then there exists a positive constant $M_{0, \nu,\nu'}$ such that:
$$
 \|\Uosc\|_{\dot{E}^s} \leq M_{0, \nu,\nu'}\Big(\|\Uoosc\|_{\dot{H}^{-1}} +\ee M_{0, \nu,\nu'}\Big)^\frac{1-s}{2} \underset{\ee \rightarrow 0}{\longrightarrow} 0,
$$
which ends the proof of the first half of Theorem \ref{Th2}.

\section{End of the proof of Theorem \ref{Th2}}

Here we prove in a direct way the convergence (in \cite{Chemin2}, $(\Ome)_{0<\ee <\ee_1}$ was proved to be a Cauchy sequence).\\

Let us emphasize that, as in Theorem \ref{Th2} the sequence of initial data $(\Uoe)_{\ee \in]0,\ee_0]}$ is assumed to be bounded in $L^2\cap \dot{H}^1$, then the same is true for its quasi-geostrophic part $\Uoqg$ (we recall that $\cQ$ is a homogeneous operator of order zero). Moreover as $\Uoqg$ goes to some quasi-geostrophic vector field $\tilde{U}_{0,QG}$ in $\dot{H}^1$ then we immediately obtain (thanks to the uniqueness of limits in the sense of distributions) that $\tilde{U}_{0,QG}$ in also in $L^2$. Next, we only have to use Theorem $2$ from \cite{FC2} claiming that System \eqref{QG1} has a unique global solution $\tilde{U}_{QG} \in \dot{E}^0 \cap \dot{E}^1$ as soon as $\tilde{U}_{0,QG} \in H^1$.

Let us consider the initial data according to the assumptions of Theorem \ref{Th2} and the unique global solution given by Theorem \ref{Th1} for a small enough $\ee$. In the previous section we already proved that the oscillating part goes to zero and we only have to study the convergence of the quasi-geostrophic part as $\ee$ goes to zero. Let us define $\dOm= \Ome-\tOm$ where $\tOm$ is the potential vorticity of the global solution $\tilde{U}_{QG}$ of the limit system. It satisfies the following system:
\begin{multline}
  \d_t \dOm +v_{\ee, QG}\cdot \n \dOm -\G \dOm= \Sum_{i=1,..., 4} B_i\\
  =-(v_{\ee, QG}-\tilde{v}_{QG})\cdot \n \tOm +(\nu-\nu') \d_3 \D \theta_{\ee, osc} +\n \Ue \cdot \n \Uosc -v_{\ee, osc}\cdot \n \Ome,
\end{multline}
supplemented by the initial data $\dOm(0)=\Omega(\Uoqg-\tilde{U}_{0,QG})$, which goes to zero in $L^2$. Taking the $L^2$ inner product with $\dOm$ we obtain:
\begin{equation}
 \frac{1}{2} \frac{d}{dt} \|\dOm\|_{L^2}^2+ \min(\nu, \nu') \|\dOm\|_{\dot{H}^1}^2 \leq \Sum_{i=1,..., 4}|(B_i, \dOm)_{L^2}|.
\label{estimdeltOm}
 \end{equation}
Three terms are easily estimated:
\begin{multline}
|(B_1, \dOm)_{L^2}|\leq \|(v_{\ee, QG}-\tilde{v}_{QG})\cdot \n \tOm\|_{\dot{H}^{-1}}\|\dOm\|_{\dot{H}^1} \leq \|v_{\ee, QG}-\tilde{v}_{QG}\|_{\dot{H}^1} \|\n \tOm\|_{\dot{H}^{-\frac{1}{2}}}\|\dOm\|_{\dot{H}^1}\\
\leq \|\dOm\|_{L^2} \|\n \tOm\|_{\dot{H}^{-\frac{1}{2}}}\|\dOm\|_{\dot{H}^1} \leq \frac{\min(\nu, \nu')}{8} \|\dOm\|_{\dot{H}^1}^2 +\frac{C}{\min(\nu, \nu')} \|\dOm\|_{L^2}^2 \|\tilde{U}_{QG}\|_{\dot{H}^\frac{3}{2}}^2.
\label{estimlim1}
\end{multline}
Similarly we get that:
\begin{equation}
\begin{cases}
\vspace{0.2cm}
  \displaystyle{|(B_3, \dOm)_{L^2}|\leq \frac{\min(\nu, \nu')}{8} \|\dOm\|_{\dot{H}^1}^2 +\frac{C}{\min(\nu, \nu')} \|\n\Ue\|_{L^2}^2 \|\n \Uosc\|_{\dot{H}^\frac{1}{2}}^2,}\\
  \displaystyle{|(B_4, \dOm)_{L^2}|\leq \frac{\min(\nu, \nu')}{8} \|\dOm\|_{\dot{H}^1}^2 +\frac{C}{\min(\nu, \nu')} \|\Ue\|_{\dot{H}^2}^2 \|\Uosc\|_{\dot{H}^\frac{1}{2}}^2.}
\end{cases}
\label{estimlim34}
\end{equation}
The last term seems more delicate at first sight, as the same problem as before appears here due to the three derivatives: we wish this term to go to zero and from the previous section $\Uosc$ goes to zero but only in $\dot{E}^s$ for $s\in [-1,1[$. As we cannot transfer more than one derivative to $\dOm$ we are stuck as $\Uosc$ is only bounded in $\dot{E}^1$. To overcome this difficulty we will simply cut the low and high frequencies in order to take advantage of both and obtain two parts that will be small each for a different reason (we refer to the appendix for the definition of the operator $\dot{S}_m$):
\begin{multline}
|(B_2, \dOm)_{L^2}|\leq |\nu-\nu'| \left(|(\d_3 \D \theta_{\ee, osc}, \dot{S}_m \dOm)_{L^2}| +|(\d_3 \D \theta_{\ee, osc}, (I_d-\dot{S}_m) \dOm)_{L^2}|\right)\\
\overset{def}{=} C_5+ C_6.
\label{C56}
\end{multline}
For the low frequencies we take advantage of the convergence of $\Uosc$ to zero:
\begin{multline}
 C_5 \leq |\nu-\nu'|\|\d_3 \D \theta_{\ee, osc}\|_{\dot{H}^{-\frac{3}{2}}}\|\dot{S}_m \dOm\|_{\dot{H}^\frac{3}{2}} \leq |\nu-\nu'|\|\Uosc\|_{\dot{H}^{\frac{3}{2}}} 2^\frac{m}{2}\|\dOm\|_{\dot{H}^1}\\
 \leq \frac{\min(\nu, \nu')}{8} \|\dOm\|_{\dot{H}^1}^2 +C 2^m\frac{|\nu-\nu'|^2}{\min(\nu, \nu')} \|\Uosc\|_{\dot{H}^\frac{3}{2}}^2.
\label{estimlim5}
 \end{multline}
Estimating the low frequencies as follows, we will not rely anymore on $\Uosc$, it is then about to choose $m$ large enough so that this term is small:
\begin{multline}
 C_6 \leq |\nu-\nu'|\|\d_3 \D \theta_{\ee, osc}\|_{\dot{H}^{-1}}\|(I_d-\dot{S}_m) \dOm\|_{\dot{H}^1}\\
 \leq |\nu-\nu'|\|\Ue\|_{\dot{H}^2} \left( \|(I_d-\dot{S}_m) \Ome\|_{\dot{H}^1} +\|(I_d-\dot{S}_m) \tOm\|_{\dot{H}^1}\right) \overset{def}{=}\|\Ue\|_{\dot{H}^2} \left(D_1+D_2\right).
\label{estimlim6}
 \end{multline}
The scheme of the proof will be, for some fixed $\eta>0$ small, to choose $m$ large enough so a part of the right-hand side from \eqref{estimdeltOm} (after time integration) is bounded by $\frac{\eta}{2}$, then to choose $\ee$ small enough so that the rest (which features in particular $2^m$ multiplied by functions going to zero as $\ee$ goes to zero) is also bounded by $\frac{\eta}{2}$. In \eqref{estimlim6}, due to $D_1$ such an $m$ \emph{a priori depends on $\ee$} which makes the previous argument impossible to perform, so we will try to cut the dependancy in $\ee$ and give a majoration by an expression going to zero as $m$ goes to infinity \emph{independantly of $\ee$}. It is not necessary for $D_2$ but for more simplicity we use the same argument for both terms, let us begin with $D_2$ as it is simpler. Thanks to the initial regularity of $\tilde{U}_{QG}$ we will not need sharp estimates, and it will be sufficient to write that (no need to introduce commutators):
\begin{equation}
 \d_t (I_d-\dot{S}_m) \tOm-\G (I_d-\dot{S}_m) \tOm =-(I_d-\dot{S}_m)\left(\tilde{v}_{QG}\cdot \n \tOm\right).
\end{equation}
Next we compute the inner product in $L^2$ with $(I_d-\dot{S}_m) \tOm$:
\begin{multline}
  \frac{1}{2} \frac{d}{dt} \|(I_d-\dot{S}_m) \tOm\|_{L^2}^2+ \min(\nu, \nu') \|(I_d-\dot{S}_m) \tOm\|_{\dot{H}^1}^2\\
  \leq \|(I_d-\dot{S}_m)\left(\tilde{v}_{QG}\cdot \n \tOm\right)\|_{\dot{H}^{-1}} \|(I_d-\dot{S}_m)\tOm\|_{\dot{H}^1}.
\end{multline}
Thanks to the classical estimates $2ab\leq a^2+b^2$ and using that for all $f\in \dot{H}^{s}\cap \dot{H}^{s+\alpha}$,
\begin{equation}
 \|(I_d-\dot{S}_m)f\|_{\dot{H}^s} \leq C \left(\int_{|\xi|\geq \frac{3}{4}2^m} |\xi|^{2s} |\hat{f}(\xi)|^2 d\xi\right)^\frac{1}{2} \leq \frac{C}{2^{\alpha m}}  \|f\|_{\dot{H}^{s+\alpha}},
 \label{majm}
\end{equation}
we get (with $s=-1$, $\alpha=\frac{1}{2}$):
\begin{multline}
\|(I_d-\dot{S}_m)\tOm (t)\|_{L^2}^2+ \min(\nu, \nu') \int_0^t\|(I_d-\dot{S}_m)\tOm (\tau)\|_{\dot{H}^1}^2 d\tau\\
\leq \|(I_d-\dot{S}_m)\tOm (0)\|_{L^2}^2 +\frac{C}{2^m}\int_0^t \|\tilde{v}_{QG}\cdot \n \tOm (\tau)\|_{\dot{H}^{-\frac{1}{2}}}^2 d\tau\\
\leq \|(I_d-\dot{S}_m)\tOm (0)\|_{L^2}^2 +\frac{C}{2^m}\int_0^t \|\tilde{v}_{QG}(\tau)\|_{\dot{H}^1}^2 \|\tilde{v}_{QG}(\tau)\|_{\dot{H}^2}^2 d\tau.
\end{multline}
Finally, thanks to the estimates provided by Theorem \ref{Th1}, we obtain:
\begin{multline}
\|(I_d-\dot{S}_m)\tOm (t)\|_{L^2}^2+ \min(\nu, \nu') \int_0^t\|(I_d-\dot{S}_m)\tOm (\tau)\|_{\dot{H}^1}^2 d\tau\\
\leq \|(I_d-\dot{S}_m)\tOm (0)\|_{L^2}^2 +\frac{C}{2^m}\frac{\|\tilde{U}_{0,QG}\|_{\dot{H}^1}^4}{\min(\nu, \nu')}.
\label{estimlim6a}
\end{multline}
We do the same for the last part (i.-e. $C_6$), starting from:
\begin{multline}
 \d_t (I_d-\dot{S}_m)\Ome-\G (I_d-\dot{S}_m)\Ome= \Sum_{i=1,2,3} E_i =\\
 -(I_d-\dot{S}_m)\left(\ve\cdot \n \Ome\right) +(\nu-\nu')\D \d_3 (I_d-\dot{S}_m)\theta_{\ee, osc}+(I_d-\dot{S}_m)\left(\n \Ue \cdot \n \Uosc\right),
\end{multline}
and as before
\begin{multline}
\|(I_d-\dot{S}_m)\Ome (t)\|_{L^2}^2+ \min(\nu, \nu') \int_0^t\|(I_d-\dot{S}_m)\Ome (\tau)\|_{\dot{H}^1}^2 d\tau\\
\leq \|(I_d-\dot{S}_m)\Ome (0)\|_{L^2}^2 +\frac{1}{\min(\nu, \nu')} \int_0^t \Sum_{i=1,2,3} \|E_i(\tau)\|_{\dot{H}^{-1}}^2 d\tau.
\end{multline}
Similarly as before (we skip details):
\begin{equation}
 \int_0^t \left(\|E_1(\tau)\|_{\dot{H}^{-1}}^2 +\|E_3(\tau)\|_{\dot{H}^{-1}}^2\right) d\tau \leq \frac{C}{2^m}\frac{\|\Uoe\|_{\dot{H}^1}^4}{\min(\nu, \nu')},
\end{equation}
and for the last term,
\begin{equation}
 \int_0^t \|E_2(\tau)\|_{\dot{H}^{-1}}^2 d\tau \leq \frac{|\nu-\nu'|^2}{\min(\nu, \nu')}\int_0^t \|(I_d-\dot{S}_m)\Ue (\tau)\|_{\dot{H}^2}^2 d\tau,
\end{equation}
we cannot perform as for the other terms as we do not have enough regularity, instead we repeat once more the same argument of truncation: applying $(I_d-\dot{S}_m)$ to \eqref{PE} we get that:
$$
\d_t (I_d-\dot{S}_m)\Ue -L (I_d-\dot{S}_m)\Ue +\frac{1}{\ee} \cA (I_d-\dot{S}_m)\Ue=\frac{1}{\ee} (-\n \Phie, 0) -(I_d-\dot{S}_m)\left(\ve\cdot \n \Ue\right),\\
$$
and computing the innerproduct in $\dot{H}^1$ with $(I_d-\dot{S}_m)\Ue$, we obtain (skipping details as they are close to the previous computations):
\begin{multline}
\|(I_d-\dot{S}_m)\Ue (t)\|_{\dot{H}^1}^2+ \min(\nu, \nu') \int_0^t\|(I_d-\dot{S}_m)\Ue (\tau)\|_{\dot{H}^2}^2 d\tau\\
\leq \|(I_d-\dot{S}_m)\Uoe\|_{\dot{H}^1}^2 +\frac{C}{\min(\nu, \nu')} \int_0^t \|(I_d-\dot{S}_m)\left(\ve\cdot \n \Ue\right)\|_{L^2}^2 d\tau\\
\leq \|(I_d-\dot{S}_m)\Uoe\|_{\dot{H}^1}^2 +\frac{C}{\min(\nu, \nu')} \frac{1}{2^m}\int_0^t \|\ve\cdot \n \Ue\|_{\dot{H}^\frac{1}{2}}^2 d\tau\\
\leq \|(I_d-\dot{S}_m)\Uoe\|_{\dot{H}^1}^2 +\frac{C}{2^m}\frac{\|\Uoe\|_{\dot{H}^1}^4}{\min(\nu, \nu')^2}.
\label{estimlim6b}
\end{multline}
Gathering \eqref{estimdeltOm}, \eqref{estimlim1}, \eqref{estimlim34}, \eqref{estimlim5} and \eqref{estimlim6} and performing an integration in time, we end up for all $t$ with:
\begin{multline}
 \|\dOm(t)\|_{L^2}^2 +\min(\nu, \nu') \int_0^t\|(\dOm(\tau)\|_{\dot{H}^1}^2 d\tau \leq \|\Uoqg-\tilde{U}_{0,QG}\|_{\dot{H}^1}^2\\
 +\frac{C}{\min(\nu, \nu')}\int_0^t \left(\|\dOm\|_{L^2}^2 \|\tilde{U}_{QG}\|_{\dot{H}^\frac{3}{2}}^2 +\|\n \Ue\|_{L^2}^2 \|\n \Uosc\|_{\dot{H}^\frac{1}{2}}^2 +\|\Ue\|_{\dot{H}^2}^2 \|\Uosc\|_{\dot{H}^\frac{1}{2}}^2 \right) d\tau\\
 +C 2^m \frac{|\nu-\nu'|}{\min(\nu, \nu')} \int_0^t \|\Uosc\|_{\dot{H}^\frac{3}{2}}^2 d\tau\\
 + |\nu-\nu'|\left(\int_0^t \|\Ue\|_{\dot{H}^2}^2 d\tau\right)^\frac{1}{2} \left( \int_0^t \Big( \|(I_d-\dot{S}_m) \Ome\|_{\dot{H}^1}^2 +\|(I_d-\dot{S}_m) \tOm\|_{\dot{H}^1}^2 \Big) d\tau \right).
\end{multline}
Thanks to the Gronwall estimates (first term in the first integral of the right-hand side), using that
$$
\int_0^t \|\tilde{U}_{QG}\|_{\dot{H}^\frac{3}{2}}^2 d\tau \leq \frac{\|\tilde{U}_{QG}\|_{L^2}\|\tilde{U}_{QG}\|_{\dot{H}^1}}{\min(\nu, \nu')},
$$
and combining it with \eqref{estimlim6a} to \eqref{estimlim6b} we obtain:
\begin{multline}
 \|\dOm(t)\|_{L^2}^2 +\min(\nu, \nu') \int_0^t\|(\dOm(\tau)\|_{\dot{H}^1}^2 d\tau \leq \exp\left(C\frac{\|\tilde{U}_{QG}\|_{L^2}\|\tilde{U}_{QG}\|_{\dot{H}^1}}{\min(\nu, \nu')^2}\right) \times \Bigg[ \|\Uoqg-\tilde{U}_{0,QG}\|_{\dot{H}^1}^2\\
 +\frac{C}{\min(\nu, \nu')^2} \|\Uosc\|_{\dot{E}^\frac{1}{2}}^2 \Big(\|\Uoe\|_{\dot{H}^1}^2+ 2^m |\nu-\nu'|^2\Big)\\
 +C\frac{|\nu-\nu'|}{\min(\nu, \nu')}\|\Uoe\|_{\dot{H}^1}\times \Big\{\|(I_d-\dot{S}_m)\tilde{\Omega}_{0, QG}\|_{L^2}^2 +\|(I_d-\dot{S}_m)\Ome (0)\|_{L^2}^2 +\frac{|\nu-\nu'|^2}{\min(\nu, \nu')}\|(I_d-\dot{S}_m)\Uoe\|_{\dot{H}^1}^2\\
 +\frac{C}{2^m} \frac{1}{\min(\nu, \nu')} \Big(\|\tilde{U}_{0,QG}\|_{\dot{H}^1}^4 +(1+ \frac{|\nu-\nu'|^2}{\min(\nu, \nu')^3}) \|\Uoe\|_{\dot{H}^1}^4\Big) \Big\}^\frac{1}{2} \Bigg].
\end{multline}
We have to be careful that in the previous estimates, for any fixed $\ee$, $\|(I_d-\dot{S}_m)\Ome (0)\|_{L^2}^2$ and $\|(I_d-\dot{S}_m)\Uoe\|_{\dot{H}^1}^2$ go to zero when $m$ goes to infinity, but nothing ensures the convergence does not depend on $\ee$. To solve the problem, we use here the extra-regularity assumption on $\Uoosc$ and \eqref{majm}:
\begin{multline}
 \|(I_d-\dot{S}_m)\Ome (0)\|_{L^2} \leq C\|(I_d-\dot{S}_m)\Uoe\|_{\dot{H}^1} \leq C\|(I_d-\dot{S}_m)\Big(\Uoosc+ (\Uoqg-\tilde{U}_{0,QG}) +\tilde{U}_{0,QG}\Big)\|_{\dot{H}^1}\\
 \leq \frac{C}{2^{\delta m}} \|\Uoosc\|_{\dot{H}^{1+\delta}} +C\|\Uoqg-\tilde{U}_{0,QG}\|_{\dot{H}^1} +C\|(I_d-\dot{S}_m)\tilde{U}_{0,QG}\|_{\dot{H}^1}\\
\leq \frac{C_0}{2^{\delta m}} +C\|\Uoqg-\tilde{U}_{0,QG}\|_{\dot{H}^1} +C\|(I_d-\dot{S}_m)\tilde{U}_{0,QG}\|_{\dot{H}^1}.
\end{multline}
To sum up we obtained that
$$
\|\dOm\|_{\dot{E}^0} \leq (1+2^m) F(\ee)+ G(m),
$$
where $F(\ee) \underset{\ee\rightarrow 0}{\longrightarrow} 0$ and $G(m) \underset{m\rightarrow \infty}{\longrightarrow} 0$. For a given $\eta>0$, let us fix $m$ large enough so that $G(m) \leq \frac{\eta}{2}$, then fix $\ee$ small enough so that $(1+2^m) F(\ee)\leq \frac{\eta}{2}$ and Theorem \ref{Th2} is proved. $\blacksquare$

\begin{rem}
 \sl{On could wonder why not using the extra regularity on $\Uoosc$ from the beginning in \eqref{C56}
 $$
 |(B_2, \dOm)_{L^2}|\leq |\nu-\nu'| \|\d_3 \D \theta_{\ee, osc}\|_{\dot{H}^{-1+\delta}}\|\dOm\|_{\dot{H}^{1-\delta}}.
 $$
 This would imply that we can prove that this additional regularity is transmitted for any time to $\Uosc$, which is not clear as taking the $\dot{H}^{1+\delta}$ inner-product of \eqref{systosc2} with $\Uosc$ we would have to deal with the term
 $$
 |\nu-\nu'|(\d_3 \Ome, \theta_{\ee, osc})_{\dot{H}^{1+\delta}}.
 $$
As $\d_3 \Ome$ can only be estimated in $L^2$, we put $2+2\delta$ derivatives on the other term which is not possible as $\theta_{\ee, osc}$ is at most in $\dot{H}^{2+\delta}$. So we are not able to estimate this term and propagate the extra regularity on $\Uosc$ unless we ask extra regularity also on the quasi-geostrophic part.
}
\end{rem}

\section{Appendix: notations}

For $s\in\R$, $\dot{H}^s$ and $H^s$ are the classical homogeneous/inhomogeneous Sobolev spaces in $\R^3$ endowed with the norms:
$$
\|u\|_{\dot{H}^s}^2=\int_{\R^3} |\xi|^{2s} |\hat{u}(\xi)| d\xi, \quad \mbox{and} \quad \|u\|_{H^s}^2=\int_{\R^3} (1+|\xi|^2)^s |\hat{u}(\xi)| d\xi.
$$
We also use the following notations: if $E$ is a Banach space and $T>0$,
$$
\cC_T{E} =\cC([0,T], E), \quad \mbox{and} \quad L_T^p{E} =L^p([0,T], E).
$$
We make abundant use of the Sobolev injections, and product laws:
\begin{prop}
 \sl{There exists a constant $C>0$ such that if $s<\frac{3}{2}$, then for any $u\in\dot{H}^s$, $u\in L^{p}(\R^3)$ with $p=\frac{6}{3-2s}$ and
 $$
 \|u\|_{L^{p}}\leq C \|u\|_{\dot{H}^s}.
 $$}
\end{prop}
\begin{prop}
 \sl{There exists a constant $C$ such that for any $(u,v)\in \dot{H}^s\times \dot{H}^t$, if $s,t<\frac{3}{2}$ and $s+t>0$ then $uv \in \dot{H}^{s+t-\frac{3}{2}}$ and we have:
 $$
 \|uv\|_{\dot{H}^{s+t-\frac{3}{2}}} \leq C \|u\|_{\dot{H}^s} \|v\|_{\dot{H}^t}.
 $$
 }
\end{prop}
Finally, we introduce the frequency truncation operator $\dot{S}^m$: consider a smooth radial function $\chi$ supported in the ball $B(0, \frac{4}{3})$, equal to $1$ in a neighborhood of $B(0, \frac{3}{4})$ and such that $r\mapsto \chi(r.e_1)$ is nonincreasing over $\R_+$. For any $u$,
$$
\dot{S}_m u =\chi(2^{-m}D) u \overset{def}{=} \cF^{-1}\Big(\chi(2^{-m}\xi) \hat{u}(\xi)\Big).
$$
This operator smoothly cuts the frequencies of size greater than $2^m$. For more details on general dyadic decompositions and Besov spaces we refer to \cite{Chemin1, Dbook}.  
\\


\textbf{Aknowledgements :} This work was supported by the ANR project INFAMIE, ANR-15-CE40-0011.

\end{document}